\documentclass[10pt]{amsart}
\usepackage{amssymb}

\usepackage{amscd}
\usepackage{amsthm}
\usepackage[dvips]{graphicx}

\usepackage[all]{xy}

\newtheorem{Thm}{Theorem}[section]
\newtheorem{Lem}[Thm]{Lemma}
\newtheorem{Def}[Thm]{Definition}
\newtheorem{Cor}[Thm]{Corollary}
\newtheorem{Prop}[Thm]{Proposition}

\newtheorem{Ex-Rem}[Thm]{Example-Remark}
\newtheorem{Rem1}[Thm]{Remark}
\newenvironment{Rem}{\begin{Rem1}\rm}{\end{Rem1}}

\title{On singular equivalences of Morita type}
\author{Guodong Zhou and Alexander Zimmermann}

\address{Guodong Zhou
 \newline Department of Mathematics
\newline East China Normal University
\newline  Dong Chuan Road 500
\newline Min Hang District
\newline Shanghai 200241
 \newline P.R.China} \email{gdzhou@math.ecnu.edu.cn}

\address{Alexander Zimmermann
\newline Universit\'e de Picardie,
\newline D\'epartement de Math\'ematiques et LAMFA (UMR 7352 du CNRS),
\newline 33 rue St Leu,
\newline F-80039 Amiens Cedex 1,
\newline France}
\email{alexander.zimmermann@u-picardie.fr}

\date{September 27, 2012}

\newcommand{\dar}{\downarrow}
\newcommand{\lra}{\longrightarrow}

\newcommand{\ra}{\rightarrow}
\newcommand{\sdp}{\times\kern-.2em\vrule height1.1ex depth-.05ex}
\newcommand{\epi}{\lra \kern-.8em\ra}

\newcommand{\N}{{\mathbb N}}

\newcommand{\Smooth}{{\mathcal{P}^{<\infty}}}

\newcommand{\ul}{\underline}
\newcommand{\ol}{\overline}

\newcommand{\Sing}{\underline{\mathrm{mod}}_{\mathcal{P}^{<\infty}}}

 \normalbaselines
\newenvironment{Proof}[1][Proof]{\begin{trivlist}
\item[\hskip \labelsep {\bfseries #1}]}{\flushright
$\Box$\end{trivlist}}

\begin{document}

\begin{abstract}
Stable equivalences of Morita type preserve  many interesting
properties and is proved to be the appropriate concept to study
for equivalences between stable categories. Recently the
singularity category attained much attraction and Xiao-Wu Chen and
Long-Gang Sun gave an appropriate definition of singular
equivalence of Morita type. We shall show that under some
conditions singular equivalences of Morita type have some
biadjoint functor properties and preserve positive degree
Hochschild homology.
\end{abstract}

\maketitle

\section*{Introduction}
For a Noetherian algebra $A$ over a commutative ring
its singularity category $D_{sg}(A)$
is defined to be
the Verdier quotient of the bounded derived category of finitely
generated modules over $A$
by the full subcategory of perfect complexes.
This notion was introduced in an unpublished
manuscript \cite{Buchweitz} by Ragnar-Olaf Buchweitz under the
name of stable derived category. He related this category to maximal
Cohen-Macaulay modules. Later Dmitri Orlov \cite{Orlov}
rediscovered this notion independently in the context of algebraic
geometry and mathematical physics, under the name
of singularity category. The derived category of an algebra is
replaced there by the derived category of coherent sheaves over a scheme.
Orlov's notation for this object seems
now to become the standard one, also in the case of the derived category of an
algebra, and we shall concentrate here on this case.

If $A$ is a selfinjective algebra, then $D_{sg}(A)$ is equivalent to the
stable category of $A$ (cf \cite{KellerVossieck, Rickardstable}).
By definition $D_{sg}(A)$ is always triangulated and
it is easy to see that $D_{sg}(A)$ is trivial if and
only if $A$ has finite global dimension.
From this point of view $D_{sg}(A)$ seems to  have advantages with
respect to the stable category of an algebra, in case the algebra
is not selfinjective, and may be an appropriate replacement.
Recently much work was undertaken to understand the structure of
$D_{sg}(A)$ under various conditions on $A$. We mention in
particular Xiao-Wu Chen's work here
\cite{ChenSchur,Chenunifying,Chenradialsquare,ChenGorenstein}, but
also Bernhard Keller, Daniel Murfet and Michel Van den
Bergh~\cite{KellerMurfetvandenBergh} as well as Osamu Iyama,
Kiriko Kato and Jun-Ichi Miyachi~\cite{IyamaKatoMiyachi}.

Abstract equivalences between stable categories of algebras are very ill-behaved, even
in case the algebras are selfinjective. Very few properties of the algebras
are preserved. However, if the equivalence is induced by an exact functor
of the module categories, much more can be said and a rich structure
is available. The concept developed for
this purpose is Brou\'e's concept of stable equivalence of Morita type~\cite{Broue}.
Since the singularity category generalises the stable category, we cannot expect better
properties in the singularity case than we have in the stable case.

Very recently analogous to the
notion of stable equivalences of Morita type, Xiao-Wu Chen and
Long-Gang Sun  defined in \cite{ChenSun} the concept of singular
equivalences of Morita type. The purpose of the present note is to study
this new concept of singular equivalences of Morita type.
We obtain two main results.
First, we shall prove in Theorem~\ref{AdjointPairs}
that under mild conditions a singular equivalence
of Morita type gives rise to a bi-adjoint pair. This section
is inspired by an analogous approach by Alex Dugas and Roberto
Martinez-Villa~\cite{DugasMartinezVilla}.
Then we shall investigate Hochschild homology and show
in Theorem~\ref{Hochschildhomology} that
Hochschild homology of a finite dimensional algebra over a field and in
strictly positive degrees
is invariant under a singular
equivalence of Morita type. The main tool here is
Serge Bouc's generalisation~\cite{Bouc}
of the Hattori-Stallings trace to Hochschild homology.

The paper is organised as follows. We recall the notion and some properties of
singularity categories in Section~\ref{singulardefinitionsection}.
Section~\ref{singularequivofmoritatype} is devoted to
the definition and some of the results of Chen and Sun on singular equivalences of
Morita type. We prove the
biadjoint property in Section~\ref{biadjointtheoremsection} and we study Hochschild
homology in Section~\ref{hochschildsection}.

\subsection*{ Acknowledgement:} We are very grateful to
Xiao-Wu Chen and Long-Gang Sun for sending us their preprint \cite{ChenSun}.

We thank the two referees for their careful reading of this paper, and
in particular for mentioning to us an error in the proof of
Theorem~\ref{AdjointPairs} which lead to a modification of the notion of being
strongly right nonsingular.

\section{Singularity  categories and singularly stable categories}
\label{singulardefinitionsection}

Let  $A$ be a right Noetherian ring. We denote by
$\mathrm{mod}(A)$ the category of finitely generated right
$A$-modules, by $D^b(\mathrm{mod}(A))$ the  bounded derived
category of $\mathrm{mod}(A)$, by $\mathcal{P}^{<\infty}(A)$ the
full subcategory of $\mathrm{mod}(A)$ consisting of modules of
finite projective dimension, and by $K^b(\mathrm{proj}(A))$ the
homotopy category of bounded complexes
of finitely generated projective $A$-modules.

\begin{Def}[\cite{Buchweitz}]
Let $A$ be a right Noetherian ring. Then the Verdier quotient
category
$$D_{sg}(A):=D^b(\mathrm{mod}(A))/K^b(\mathrm{proj}(A))$$ is called the
singularity category of $A$.
\end{Def}

It is well-known that $K^b(\mathrm{proj}(A))$ is a full triangulated
subcategory of $D^b(A)$.  We briefly recall the construction
of the Verdier quotient.
We refer to Gabriel and Zisman's book \cite[Chapter 1]{GabrielZisman} for
more ample details, and give only the basic construction here for the
convenience of the reader.

The objects of $D_{sg}(A)$ are the same as those
of $D^b(A)$. Let $X$ and $Y$ be objects of $D_{sg}(A)$. Then a
morphism in $Hom_{D_{sg}(A)}(X,Y)$ is represented by triples
$(\nu,Z,\alpha)$ where $Z$ is an object in $D^b(A)$, where
$\alpha\in Hom_{D^b(A)}(Z,Y)$ and where $\nu\in Hom_{D^b(A)}(Z,X)$
so that the mapping cone of $\nu$ is isomorphic to an object in
$K^b(\mathrm{proj}(A))$. A triple $(\nu,Z,\alpha)$ is covered by a
triple $(\nu',Z',\alpha')$ if there is a morphism $\psi\in
Hom_{D^b(A)}(Z',Z)$ so that $\nu'=\nu\circ\psi$ and
$\alpha'=\alpha\circ\psi$. Two triples $(\nu,Z,\alpha)$ and
$(\nu'',Z'',\alpha'')$ are equivalent if both are covered by some
triple $(\nu',Z',\alpha')$. This way the category of triples is
directed, and the morphisms from $X$ to $Y$ is the limit of this
category.

The construction of the singularity category as Verdier quotient
implies that $D_{sg}(A)$ is always triangulated.

Let $A$ be any right Noetherian ring. Denote by
$\ul{\mathrm{mod}}(A)$ the stable category of (finitely generated
right) $A$-modules, with objects being the same as
$\mathrm{mod}(A)$ and morphisms $\ul{\mathrm{Hom}}_A(M,N)$ being
the equivalence classes of morphisms of $A$-modules modulo those
factoring through a projective module. Recall that the category
$\ul{\mathrm{mod}}(A)$ admits an endo-functor $\Omega$, the syzygy
functor, defined as $\ker(\pi_X)$, where for every object $X$ in
$\ul{\mathrm{mod}}(A)$ we choose a projective object $P_X$ in
$\mathrm{mod}(A)$ and an epimorphism $P_X\stackrel{\pi_X}{\lra}X$
in $\mathrm{mod}(A)$.

By the very construction there are natural functors
\begin{eqnarray*}\mathrm{mod}(A)&\stackrel{F}{\lra}& \ul{\mathrm{mod}}(A)\\
D^b(A)&\stackrel{G}{\lra}& D_{sg}(A)\\
\ul{\mathrm{mod}}(A)&\stackrel{H}{\lra}& D_{sg}(A)
\end{eqnarray*}
so that the diagram
$$\begin{array}{ccc}
\mathrm{mod}(A)&\stackrel{}{\lra}&D^b(A)\\
\dar F&&\dar G\\
\ul{\mathrm{mod}}(A)&\stackrel{H}{\lra}& D_{sg}(A)
\end{array}$$
commutes.
Moreover,
$H(M)=0$ if and only if $M$ is of finite projective dimension. Finally
$H$ commutes with syzygies
in the sense that
$$ H\circ\Omega\simeq [-1]\circ H.$$
A consequence of this relation is an important  observation,
namely that the singularity category is in general not Hom-finite
\cite{Chenradialsquare} and is in general not a Krull-Schmidt
category. An example is given by the $3$-dimensional local algebra
$A=K[X,Y]/(X^2,Y^2,XY)$ over a field $K$. Indeed, for the simple
$A$-module $S$ one gets $\Omega(S)\simeq S\oplus S$ and this
implies isomorphisms
$$H(S)\simeq H(\Omega(S))[1]\simeq H(S)[1]\oplus H(S)[1]\simeq
(H(S)[2])^4\simeq \cdots \simeq (H(S)[n])^{2^n}$$
in $D_{sg}(A)$. Moreover, $H(S)\neq 0$ since $S$ is of infinite
projective dimension. Therefore,
$$\mathrm{dim}_K\mathrm{End}_{D_{sg}(A)}(H(S))=
\mathrm{dim}_K\mathrm{End}_{D_{sg}(A)}( (H(S)[n])^{2^n})\geq 2^n$$ for all
$n\geq 0$ and this implies that $\mathrm{dim}_K
\mathrm{End}_{D_{sg}(A)}(H(S))=+\infty$.

As we have seen, $M$ is an $A$-module of finite projective
dimension if and only if $H(M)=0$. Hence, it is natural to consider the
following category $\Sing(A)$.

\begin{Def} Let $A$ be a finite dimensional algebra. The
singularly
stable category is by definition the quotient category of
$\mathrm{mod}(A)$ by $\mathcal{P}^{<\infty}(A)$, denoted by
$\Sing(A):= \mathrm{mod}(A)/\mathcal{P}^{<\infty}(A)$.

More precisely, the objects of $\Sing(A)$ are the same objects as
those in $\mathrm{mod}(A)$ and for two $A$-modules $X$ and $Y$
define $\mathrm{Hom}_{\Sing(A)}(X,Y)$ to be the equivalence
classes of $A$-module homomorphisms $X\lra Y$ modulo those
factoring through an object in $\mathcal{P}^{<\infty}(A)$.
\end{Def}

It is clear that $H$ factors through the natural functor
$$\ul{\mathrm{mod}}(A)\stackrel{\Pi}{\lra} \Sing(A)$$
in the sense that there is a natural functor
$$\Sing(A)\stackrel{L}{\lra}D_{sg}(A)$$
so that $$H=L\circ\Pi.$$

\begin{Rem}\label{Remexample}
Observe that $L$ is not an embedding in general.
Let $Q$ be the quiver \unitlength=1cm
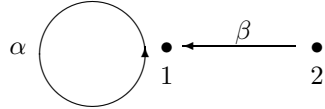
\begin{center}\begin{picture}(5,1.5)
\put(2,.6){$\bullet$}\put(4,.6){$\bullet$}
\put(2,.2){$1$}\put(4,.2){$2$} \put(3.8,.7){\vector(-1,0){1.5}}
\put(1.1,.6){\circle{2.4}} \put(1.8,.6){\vector(0,1){.1}}
\put(0,.6){$\alpha$} \put(3,.8){$\beta$}
\end{picture}\end{center}
and let $A=KQ/\langle \alpha^2,\beta\alpha\rangle$. Let $S_1$ and
$S_2$ be the two simple $A$-modules. Then $H(S_1)\simeq H(S_2)$
since $\Omega^2(S_2)\simeq S_1\simeq  \Omega^2(S_1)$,
but $S_1\not\simeq S_2$ in $\Sing(A)$ since there is no
non zero homomorphism of $A$-modules between these objects.
\end{Rem}

\begin{Rem}
We could consider modules of finite projective dimension as
``smooth" objects. Then the singularly stable category measures the
singularity of $A$.  Clearly the algebra $A$ has finite global
dimension if and only if the singularly stable category has only
one object with only one endomorphism. However, the singularly
stable category is only an additive category, and in general it is
not triangulated.
If $A$ is selfinjective, $H$ is an equivalence
(cf \cite{KellerVossieck,Rickardstable})
and an $A$-module of finite projective dimension is actually projective.
Hence also $\Pi$ is an equivalence in this case.
\end{Rem}

\begin{Rem}
Let $A$ be the algebra introduced in Remark~\ref{Remexample}. Then
it is easy to see that $D_{sg}(A)\simeq D_{sg}(K[X]/X^2)$.
However, these two algebras are not stably equivalent.
In fact, if  they were stably equivalent, then there would be
a one to one correspondence
between the isomorphism classes of non-projective indecomposable modules.
However, up to
isomorphisms, $A$ has more than two non-projective indecomposable modules
and $K[X]/(X^2)$ has only one such module.

We are grateful to one of the referees who suggested the above proof
which is simpler than our original argument.
\end{Rem}

\section{Singular equivalences of Morita type}
\label{singularequivofmoritatype}

As mentioned in the introduction, general stable
equivalences have very poor properties, even for selfinjective algebras.
A richer concept is given by Brou\'e~\cite{Broue}.
Brou\'e defined stable equivalences of Morita type
as equivalences between stable module categories
induced by tensor product with bimodules. This concept was highly successful
in the understanding of equivalence between stable categories
of self-injective algebras and was a subject of numerous studies.

We consider the question when the singularly stable categories of
two algebras are equivalent. Since equivalences between singular categories
of selfinjective algebras coincide with stable equivalences, we need a
richer concept than just an equivalence between triangulated categories.
Recently Xiao-Wu Chen and Long-Gang
Sun introduced singular equivalences of Morita type \cite{ChenSun}
on the model of Brou\'e's concept of stable equivalences of Morita
type.

Let $K$ be a commutative ring. For a
$K$-algebra $A$, we denote by $A^e=A^{op}\otimes_K A$ its
enveloping algebra.

\begin{Def} (cf \cite{ChenSun}) \label{singularequivalenceofMoritatypedef}
Let $A$ and $B$ be two $K$-algebras for a commutative ring $K$.
Let $_AM_B$ and $_BN_A$ be two bimodules so that
\begin{itemize}
\item $M$ is finitely generated and projective as $A^{op}$-module
and as $B$-module; \item $N$ is finitely generated and projective
as $A$-module and as $B^{op}$-module; \item
${}_AM\otimes_BN_A\simeq {}_AA_A\oplus {}_AX_A$ for a module $X\in
\mathcal{P}^{<\infty}(A^e)$; \item ${}_BN\otimes_AM_B\simeq
{}_BB_B\oplus {}_BY_B$ for a module $Y\in
\mathcal{P}^{<\infty}(B^e)$.
\end{itemize}
We then say that the pair  $({}_AM_B, {}_BN_A)$
induces a singular equivalence of Morita type.

We say that $A$ and $B$ are singularly equivalent
of Morita type if there is a pair of bimodules $({}_AM_B, {}_BN_A)$
which induces a singular equivalence of Morita type.
\end{Def}

\begin{Rem}\begin{itemize}
\item
It is immediate from the definition that a pair of
bimodules inducing a stable equivalence of Morita type induces a singular
equivalence of Morita type as well.

However, a singular equivalence of Morita type will not be a
stable equivalence of Morita type in general since the property of
$X$ to be in $\Smooth (A^e)$ is in general much weaker than the
condition to be projective as bimodule.

Nevertheless, if $A$ is
selfinjective (and thus so is any algebra singularly equivalent of Morita type
to $A$,  as is remarked in \cite{ChenSun}),
any module with finite projective resolution is actually projective,
and hence
a singular equivalence of Morita type is actually a stable equivalence
of Morita type. The concept of a singular equivalence of Morita type and
of a stable equivalence of Morita type coincide for selfinjective algebras.
\item
Let $({}_AM_B, {}_BN_A)$ induce a singular equivalence of Morita type and let
$M\otimes_BN\simeq A\oplus X$ and $N\otimes_AM\simeq B\oplus Y$.
Then $X$
is projective as $A$-left module and as $A$-right module. Indeed,
$M$ is projective as $B$-right module, hence a direct factor of some $B^n$.
Hence $M\otimes_BN$ is a direct factor of $B^n\otimes_BN\simeq N^n$.
Now, $X$ is by definition a direct factor of $N^n$ and
since $N$ is projective as $A$-right module, $X$
is projective as $A$-right module.
Similarly $X$ is projective on the left.
Likewise $Y$ is projective as $B$-left module and as $B$-right module.
\end{itemize}
\end{Rem}

From now on to the end of the present section and in
Section~\ref{biadjointtheoremsection} fix a field $K$
and $K$-algebras will always be
supposed to be finite dimensional and modules will be always
finitely generated.

The following result is a direct consequence of
Definition~\ref{singularequivalenceofMoritatypedef}.

\begin{Prop}
Let  $({}_AM_B,  {}_BN_A)$ be a pair of bimodules inducing a
singular equivalence of
Morita type between two $K$-algebras $A$ and $B$. Then
$$-\otimes_AM_B  : D_{sg}(A)\lra D_{sg}(B)$$
is an equivalence of
triangulated categories  with quasi-inverse
$$-\otimes_B N_A:D_{sg}(B) \lra D_{sg}(A).$$
Moreover, the same functors establish
an equivalence of additive   categories between
$\underline{\mathrm{mod}}_{\mathcal{P}^{<\infty}}(A)$ and
$\underline{\mathrm{mod}}_{\mathcal{P}^{<\infty}} (B)$.
\end{Prop}

The following result is an adaptation to the singular situation of a
proof of Yu-Ming Liu for stable equivalences of Morita type (cf
\cite[Lemma 2.2]{Liu2004}). The proof carries over verbatim.

\begin{Prop}\label{progenerator}(cf \cite{ChenSun})
Let $A$ and $B$ be $K$-algebras. Suppose $({}_AM_B,  {}_BN_A)$
induces a singular equivalence of Morita type. Then $_AM$  is a
progenerator in $mod(A^{op})$, and likewise for $M_B$, ${}_BN$ and
$N_A$. \end{Prop}

The following fact is proved in \cite{ChenSun} analogous to
\cite[Proposition 2.1 and  Theorem 2.2]{Liu2008}.

\begin{Prop}(cf \cite{ChenSun})\label{blocwiseequivalence}
Let $A$ and $B$ be two $K$-algebras without direct summands which
have finite projective
dimension as bimodules. Assume that two bimodules ${}_AM_B$ and
${}_BN_A$ induce a singular equivalence of Morita type between $A$
and $B$.

\begin{enumerate}
\item Then $A$ and $B$ have the same number of indecomposable
summands. In particular, $A$ is indecomposable if and only if $B$
is indecomposable.

\item  Suppose that $A = A_1 \times A_2\times \cdots \times
A_s$ and $B = B_1\times B_2\times \cdots\times B_s$, where all
$A_i$ and all $B_i$ are indecomposable algebras. Then, there is a permutation
$\sigma$ of the set $\{1,\dots,s\}$ so that
$A_i$ and $B_{\sigma(i)}$ are singularly
equivalent of Morita type for all $i\in\{1,\dots,s\}$.
\end{enumerate}
\end{Prop}

In analogy of what is known to hold for
stable equivalences of Morita type,
Chen and Sun also show the following lemma.

\begin{Lem}(cf \cite{ChenSun}) \label{uniquedirectfactor}
Let $K$ be a field and let
$A$ and $B$ be finite dimensional $K$-algebras.
Assume that bimodules ${}_AM_B$ and ${}_BN_A$ define a singular
equivalence of Morita type between $A$ and $B$, and suppose
that $A$ or $B$ is indecomposable as an algebra. Then $M$ and $N$
each have a unique indecomposable bimodule summand of infinite
projective dimension. If we denote these summands as $M_1$ and
$N_1$ respectively, then $(M_1,N_1)$ also induces a
singular equivalence of Morita type between $A$ and $B$.
\end{Lem}

Let $K$ be a field, and let $A$ and $B$ be finite dimensional
$K$-algebras without direct summands having   finite projective
dimension as bimodules. Proposition~\ref{blocwiseequivalence} and Lemma~\ref{uniquedirectfactor} imply that for a singular
equivalence of Morita
type induced by $({}_AM_B, {}_BN_A)$ we can always suppose
that $A$ and $B$ are
indecomposable algebras and that ${}_AM_B$ and  ${}_BN_A$ are
indecomposable bimodules.

\begin{Rem} During the ICRA 2012 in Bielefeld,
Chang-Chang Xi raised the question whether there are algebras which are
singularly equivalent of Morita type, but  which are {\em not}
stably equivalent of Morita type. This remark answers this question.

For any algebra $A$ denote by
$T_2(A):=\left(\begin{array}{cc} A & A\\ 0& A\end{array}\right)$
the algebra of upper $2\times 2$ triangular matrices over $A$.
In a forthcoming paper, Yu-Ming Liu and the first author give two
indecomposable $K$-algebras $A$ and $B$ which are stably equivalent
but not Morita equivalent, but for which  $T_2(A)$
is not stably equivalent to  $T_2(B)$.
However, Chen and Sun (\cite{ChenSun}) show that if $A$ and
$B$ are singular equivalent of Morita type, then also
$T_2(A)$ and $T_2(B)$ are singular equivalent of Morita type.
\end{Rem}

\section{Singular equivalences of Morita type give adjoint pairs}

\label{biadjointtheoremsection}

Our aim is to prove analogous result of
Dugas and Martinez-Villa \cite[Theorem 2.7]{DugasMartinezVilla}
for singular equivalences of Morita type. For a $K$-algebra $A$,
denote by $J(A)$ its Jacobson radical.

\begin{Thm}\label{AdjointPairs}
Let $K$ be a field and let $A$ and $B$ be finite dimensional
indecomposable $K$-algebras. Suppose $A$ and $B$ are not of finite projective
dimension as bimodules and suppose that $A/J(A)$ and $B/J(B)$ are
separable over $K$.  Let  $({}_AM_B, {} _BN_A)$ be a pair of bimodules
inducing a singular equivalence of Morita type between $A$ and $B$.
Suppose that  ${}_AM_B$ is indecomposable as a bimodule, and suppose that
$Hom_{A^{op}}({}_AM_B,{}_AA_A)$ is projective as a $B^{op}$-module.

Then
$${}_BN_A\simeq Hom_{A^{op}}({}_AM_B, {}_AA_A)$$
as $B^{op}\otimes_KA$-modules, and $(-\otimes_BN, -\otimes_AM)$ is
a pair of adjoint functors between the module categories $\mathrm{mod}(B)$
 and  $\mathrm{mod}(A)$.
\end{Thm}

\begin{Rem}
Since a singular equivalence of Morita type induces an equivalence
$D_{sg}(A)\simeq D_{sg}(B)$ and $\Sing(A)\simeq \Sing(B)$
it is clear that $(-\otimes_AM,-\otimes_BN)$ is a pair of adjoint functors
between $D_{sg}(A)$ and $D_{sg}(B)$, as well as between
$\Sing(A)$ and $\Sing(B)$. Theorem~\ref{AdjointPairs} states
that the functors form an
adjoint pair {\em between the module categories}.
\end{Rem}

In order to prove Theorem~\ref{AdjointPairs}, we shall use
the following technical notion, motivated by Dugas and
Martinez-Villa \cite{DugasMartinezVilla}.

\begin{Def}
An $A^{op}\otimes_KB$-module ${}_AU_B$
is called  strongly  right nonsingular, if for each
$A$-module $T_A$,  the $B$-module $\Omega^n_A(T)\otimes_AU_B$ is projective for $n>>0$.

\end{Def}

\begin{Lem}\label{RSN}
Let $K$ be a field and let $A$ and $B$ be finite dimensional $K$-algebras.
\begin{itemize}

\item[(i)] Let ${}_AU_B$ be a bimodule which is projective as a left
and as a right module.  Then ${}_AU_B$ is  strongly  right nonsingular
if and only if for each $A$-module $T_A$,
 $T\otimes_A U_B$ has finite  projective dimension.

\item[(ii)] Objects in $\Smooth(A^{op}\otimes_K B)$ which are
projective as a left and as  a right modules are   strongly right nonsingular. In particular, for a singular equivalence of Morita
type induced by the pair of bimodules $({}_AM_B,{}_BN_A)$, so that
$M\otimes_BN\simeq A\oplus X$ in $\mathrm{mod}(A\otimes_KA^{op})$
and $N\otimes_AM\simeq B\oplus Y$ in
$\mathrm{mod}(B\otimes_KB^{op})$, the two bimodules ${}_AX_A$ and
${}_BY_B$ are   strongly right nonsingular bimodules.

\item[(iii)] Let $A$ be an algebra such that $A/J(A)$ is separable
over $K$. If the $A^e$-module $A$ is not in $\Smooth(A^e)$,
then the bimodule ${}_AA_A$ is not   strongly right nonsingular.

\item[(iv)] Let $A$ and $B$ be finite dimensional indecomposable
$K$-algebras which are not of finite projective dimension as
bimodules and such that $A/J(A)$ and $B/J(B)$ are separable over
$K$. Let $({}_AM_B, {} _BN_A)$ be a pair of bimodules
inducing a singular equivalence of
Morita type between $A$ and $B$.   Then the two bimodules ${}_AM_B$
and ${}_BN_A$ are not  strongly right nonsingular.

\item[(v)] A direct summand of a  strongly  right nonsingular
bimodule is also  strongly  right nonsingular.
The direct sum
of two right strongly non singular bimodules is also
strongly  right nonsingular.

\end{itemize}
\end{Lem}

\begin{Proof}

(i). Suppose that  for each $A$-module $T_A$,
$T\otimes_A U_B$ has finite projective
dimension. Then for an
$A$-module $T_A$, take a minimal  projective resolution
$$\cdots \to P_n\to P_{n-1}\to \cdots \to P_1\to P_0\to T_A\to 0$$
and apply $-\otimes_AU_B$. The result is a complex of $B$-modules
$$\cdots \to P_n\otimes_A U_B\to P_{n-1}\otimes_A U_B\to
\cdots \to P_1\otimes_A U_B\to P_0\otimes_A U_B\to T\otimes_A U_B
\to 0. $$ This complex is actually exact, as ${}_AU$ is
projective. For $n\geq 1$, we have an exact sequence
$$0\to \Omega^n_A(T_A) \otimes_A U_B\to P_{n-1}\otimes_A U_B\to
\cdots \to P_1\otimes_A U_B\to P_0\otimes_A U_B\to T\otimes_A U_B
\to 0. $$
Note that for $0\leq i\leq n-1$, $P_i\otimes_AU_B$ is
projective, as $U_B$ is projective. Since
$T\otimes_A U_B$ has finite projective dimension, by
Schanuel's Lemma,  for  $n>> 0$ we get that
$\Omega^n_A(T)\otimes_AU_B$ is projective as a
$B$-module. This proves that ${}_AU_B$ is strongly right nonsingular.

Conversely, suppose that ${}_AU_B$ is strongly right nonsingular.
Take a minimal projective resolution
$$\cdots \to P_n\to P_{n-1}\to \cdots \to P_1\to P_0\to T_A\to 0$$
and apply $-\otimes_AU_B$ to get a complex
$$0\to \Omega^n_A(T)\otimes_A U_B\to P_{n-1}\otimes_A U_B\to
\cdots \to P_1\otimes_A U_B\to P_0\otimes_A U_B\to T\otimes_A U_B
\to 0 $$ of $B$-modules. This complex is exact, as ${}_AU$ is
projective.   As for $n>>0$, we have that $\Omega^n_A(T)\otimes_A U_B$ is projective,
$T\otimes_A U_B$ has finite projective dimension.

We shall use (i) in the proof of (ii)-(iv).
\medskip

(ii). Let ${}_AU_B$ be a bimodule of finite projective dimension
which is projective as left and as right module. Then there exists an
exact sequence of $A^{op}\otimes_KB$-modules
$$0\to P_n\to P_{n-1}\to \cdots \to P_0\to U\to 0,$$
where for any $0\leq i\leq n$, $P_i$ is a projective
$A^{op}\otimes_KB$-module. As ${}_AU$ is projective, the above
sequence splits  as exact sequence of left modules. So if we apply
$T_A\otimes_A-$, it remains exact. Observe that all the
$B$-modules $(T\otimes_A P_i)_B$ are projective and thus the
$B$-module $(T\otimes_A U)_B$ has finite projective dimension. We
have proved that ${}_AU_B$ is  strongly right  nonsingular.  The
second statement follows from the first one by observing that the
two bimodules ${}_AX_A$ and ${}_BY_B$ are projective as left and
as right modules.

\medskip

(iii).
Suppose ${}_AA_A$ is
 strongly right nonsingular. Then by (i) for each right $A$-module $T_A$,
the module $T_A\simeq T\otimes_AA_A$ is of finite projective
dimension. Therefore $A$ has finite global dimension and by
\cite[Section 1]{BGMS}, we have  $A\in \Smooth(A^e)$.
This proves the statement. Note that the relevant conclusion from
\cite[Section 1]{BGMS} is shown only under the hypothesis that
$A/J(A)$ is separable.

\medskip

(iv). For
each right $B$-module $T_B$ we get isomorphisms of $B$-modules
$$T\otimes_B N\otimes_A M_B\simeq
T\otimes_B({}_BB_B\oplus {}_BY_B)\simeq T_B\oplus (T\otimes_B
Y_B).$$ If ${}_AM_B$ is  strongly right nonsingular, by (i)
$T\otimes_B N\otimes_A M_B$
has finite projective dimension as a right $B$-module, and thus $T_B$ has finite
projective dimension. As in (iii), this implies that the
$B^e$-module $B$ is an object of $\Smooth(B^e)$, which is a
contradiction to the hypothesis on $B$.

The case of ${}_BN_A$ is similar.

\medskip

(v) is trivial.

\end{Proof}

\begin{Rem}
In \cite[Section 1]{BGMS} an example is given showing that we
do need the hypothesis in (ii) and (iii) that $A/J(A)$ and $B/J(B)$
are separable over $K$.
\end{Rem}

\begin{Proof}[Proof of Theorem~\ref{AdjointPairs}]

Denote ${}_B\check M_A:=\mathrm{Hom}_{A^{op}}({}_AM_B,{}_AA_A)$
to simplify the notation.
Then $(-\otimes_B \check M_A,-\otimes_AM_B)$ is an adjoint pair of
functors between $\mathrm{mod}(B)$ and $\mathrm{mod}(A)$, because
$$\mathrm{Hom}_A({}_B\check M_A, {}_AA_A)=\mathrm{Hom}_A(\mathrm{Hom}_A({}_AM_B,
{}_AA_A), {}_AA_A)\simeq {}_AM_B,$$ as ${}_AM$ is finitely
generated projective.  This pair of adjoint functors can be
defined on $D^b(\mathrm{mod}(B))$ and $D^b(\mathrm{mod}(A))$, as
${}_B\check M$ and ${}_AM$ are finitely generated projective
modules. They induce functors between $D_{sg}(A)$ and $D_{sg}(B)$
since $-\otimes_AM_B$ maps $K^b(\mathrm{proj}(A))$ to
$K^b(\mathrm{proj}(B))$, and since $-\otimes_B\check M_A$ maps
$K^b(\mathrm{proj}(B))$ to $K^b(\mathrm{proj}(A))$.

Moreover, $-\otimes_AM_B$ and $-\otimes_B\check M_A$
induce functors between $\Sing(A)$ and $\Sing(B)$
since $\Smooth(A)\otimes_AM_B$ belongs to $\Smooth(B)$ and
likewise for $-\otimes_B\check M_A$.

Let $\eta:id_{\mathrm{mod}(B)}\lra -\otimes_B\check M\otimes_A M_B$
be the unit of the
adjoint pair $(-\otimes_B \check M_A,-\otimes_AM_B)$ between
$\mathrm{mod}(B)$ and $\mathrm{mod}(A)$ and let
$\eta_B: B\to {}_B\check M\otimes_A M_B$ be its
evaluation on $B$.
As ${}_B\check
M\otimes_A M_B\simeq \mathrm{End}_A({}_AM_B)$ as $B^e$-modules,
$\eta_B$ identifies with the structure map of the right $B$-module
structure of $M$. By Lemma~\ref{progenerator}, $\eta_B$ is
injective and we can form a short  exact sequence as follows:
$$0\to {}_BB_B\stackrel{\eta_B}{\to} {}_B\check M\otimes_A M_B\to
{}_BU_B\to 0\ \ (*).$$

Applying ${}_AM_B\otimes_B-$ to the exact sequence (*)
gives the exact sequence
$$0\to {}_AM_B\stackrel{Id_M\otimes \eta_B}{\to}
{}_AM\otimes_B\check M\otimes_A M_B\to {}_AM\otimes_BU_B\to 0.$$
Now it is easy to see that the monomorphism
$Id_M\otimes \eta_B$ is split by the bimodule map
$${}_AM\otimes_B\check M\otimes_AM\simeq
{}_AM\otimes_B\mathrm{End}_A({}_AM_B)\to {}_AM_B$$
where the
second map is the evaluation map. Hence
$${}_AM\otimes_B\check M\otimes_AM_B\simeq {}_AM_B\oplus ({}_AM\otimes_BU_B).$$

\noindent
\textbf{Claim 1}:  $U_B$ is projective and  ${}_BU_B$ is
strongly right nonsingular.

We shall use this claim for the moment and  and give the proof of
Claim 1 just after having finished the proof of
Theorem~\ref{AdjointPairs}.

\medskip

Applying $-\otimes_BN$ to the isomorphism
$${}_AM\otimes_B\check M\otimes_AM_B\simeq {}_AM_B\oplus ({}_AM\otimes_BU_B)$$
gives
$$\begin{array}{rcl} {}_AM\otimes_B\check
M\otimes_AM\otimes_BN_A &\simeq &({}_AM_B\oplus ({}_AM\otimes_BU_B))\otimes_BN_A\\
&\simeq & ({}_AM\otimes_BN_A) \oplus
({}_AM\otimes_BU\otimes_BN_A)\\
&\simeq &  {}_AA_A\oplus {}_AX_A\oplus
({}_AM\otimes_BU\otimes_BN_A). \end{array}$$
But we also get

$$\begin{array}{rcl}
{}_AM\otimes_B\check
M\otimes_AM\otimes_BN_A &\simeq &{}_AM\otimes_B\check M\otimes_A(M\otimes_BN_A)\\
&\simeq &{}_AM\otimes_B\check
M\otimes_A({}_AA_A\oplus {}_AX_A) \\
&\simeq & ({}_AM\otimes_B\check M_A) \oplus ({}_AM\otimes_B\check
M\otimes_A X_A).
\end{array}$$

\noindent
\textbf{Claim 2}:   ${}_AM\otimes_BU\otimes_BN_A$ and
${}_AM\otimes_B\check M\otimes_A X_A$ are  strongly right
nonsingular.

Again we shall use this claim for the moment and give the proof
of Claim 2 just after having finished the proof of Theorem~\ref{AdjointPairs}.

\medskip

The indecomposable $A^e$-module $A$ is not  strongly right  nonsingular by Lemma~\ref{RSN} part (iii). The Krull-Schmidt theorem
shows that the $A^e$-module $A$ is a direct factor of
$M\otimes_B\check M$ or of ${}_AM\otimes_B\check M\otimes_A X_A$.
Claim 2  shows that ${}_AM\otimes_B\check M\otimes_A X_A$ is
strongly  right nonsingular, and hence all of its direct factors. Hence
the $A^e$-module $A$ is a direct factor of $M\otimes_B\check M$.
This shows that there is an $A^e$-module $\tilde X$ such that
$${}_AM\otimes_B\check
M_A\simeq {}_AA_A\oplus {}_A\tilde{X}_A\ \ (**).$$
The
bimodule ${}_A\tilde{X}_A$ is  strongly  right nonsingular by
Lemma~\ref{RSN} (ii) and (v), as ${}_A\tilde{X}_A$ is a direct
summand of ${}_AX_A\oplus ({}_AM\otimes_BU\otimes_BN_A)$.

Now we apply $N\otimes_A-$ to (**)  and get
$${}_BN\otimes_A M\otimes_B\check M_A\simeq {}_BN_A \oplus
({}_BN\otimes_A\tilde{X}_A),$$
but
$${}_BN\otimes_A M\otimes_B\check M_A
\simeq
({}_BB_B\oplus {}_BY_B)\otimes_B\check M_A\simeq
{}_B\check M_A\oplus ({}_BY\otimes_B\check M_A).$$
So
$${}_BN_A
\oplus ({}_BN\otimes_A\tilde{X}_A)\simeq {}_B\check M_A\oplus
({}_BY\otimes_B\check M_A).$$

\noindent
\textbf{Claim 3}: ${}_BN\otimes_A\tilde{X}_A$ and
${}_BY\otimes_B\check M_A$ are  strongly  right nonsingular; the $B^{op}\otimes_KA$-module
${}_B\check M_A$ is indecomposable.

Again we shall use this claim for the moment and give the proof
of Claim 3 just after having finished the proof of Theorem~\ref{AdjointPairs}.

\medskip

As in Lemma~\ref{RSN} (iv) the module  ${}_BN_A$ is not  strongly right
nonsingular. We hence
obtain that the two indecomposable bimodules ${}_BN_A$ and
${}_B\check M_A$ are isomorphic.
\end{Proof}

\begin{Proof}[Proof of Claim 1]  As in the paragraph preceding the statement
of Claim 1,  we have an isomorphism of bimodules
$${}_AM\otimes_B\check M\otimes_AM_B\simeq {}_AM_B\oplus ({}_AM\otimes_BU_B).$$
Since $M_B$ and $\check M_A$ are projective,  $M\otimes_BU_B$ is projective
as a  right $B$-module and since $M_B$ is a progenerator by
Proposition~\ref{progenerator}, we see that $U_B$ is projective.

Given a right $B$-module $T_B$,
we apply $T\otimes_B-$ to (*)  and we get an exact sequence
$$T_B\stackrel {\eta_T}{\to} T\otimes_B \check
M\otimes_A M_B\to T\otimes_B U_B\to 0,$$ where
$\eta_T=id_T\otimes_B\eta_B$.

As $\eta_T$ is an isomorphism in
$D_{sg}(B)$, there exists $n>>0$ such that $\Omega^n(\eta_T)$ is
an isomorphism in $\underline{\mathrm{mod}}(B)$.
In fact, by \cite[Example
2.3]{KellerVossieck} or \cite[Corollary 3.9(1)]{Beligiannis},
given two $B$-modules $V$ and $W$, we have
$$\mathrm{Hom}_{D_{sg}(B)}(V, W)=\underrightarrow{\mathrm{lim}}\
\ul{\mathrm{Hom}}_B (\Omega^iV, \Omega^i W).$$
Suppose that a module homomorphism $f:V\to W$ is invertible
in the singularity category.  Then its inverse  is
induced from a module homomorphism $g: \Omega^i(W)\to \Omega^i(V)$.
We see that $\Omega^i(f)\circ g$ coincides with $Id_W$
(resp. $g\circ \Omega^i(f) $ coincides with $Id_V$) in the singularity
category, so   $\Omega^{n-i}(\Omega^i(f)\circ g)=\Omega^n(f)\circ \Omega^{n-i}(g)$
coincides with $Id_{\Omega^n(N)}$  in
$\ul{\mathrm{Hom}}_B(\Omega^n(V), \Omega^n(W))$ for big enough $n$.

Let $P_*$ be the minimal projective resolution of $T_B$ and let $Q_*$ be
the minimal projective resolution of $T\otimes_B\check
M\otimes_AM_B$. As $P_*\otimes_B\check M\otimes_AM_B$ is also a
projective resolution of $T\otimes_B\check M\otimes_AM_B$, the
Comparison Lemma gives a chain map $f_*:P_*\otimes_B\check
M\otimes_AM_B\to Q_*$. Therefore, we have a commutative diagram
$$\xymatrix{ P_*\ar[r] \ar[d]_{\eta_{P_*}}& T_B\ar[d]_{\eta_T}\\ P_*\otimes_B\check
M\otimes_AM_B\ar[r] \ar[d]_{f_*} &T\otimes_B\check M\otimes_AM_B\ar[d]_{=}\\
Q_*\ar[r] &T\otimes_B\check M\otimes_AM_B}$$

Note that the induced map
$$\Omega^n_B(T_B)\stackrel{\eta_{\Omega^n_B(T_B)}}{\to}
\Omega^n_B(T_B)\otimes_B\check
M\otimes_AM_B\stackrel{f_n}{\to} \Omega^n_B(T\otimes_B\check
M\otimes_AM_B)$$ is just $\Omega^n(\eta_T)$, which is an isomorphism
as $n$ is supposed to be large enough, as we have seen.

As $f_n$ induces an isomorphism between
$\Omega^n_B(T\otimes_B\check M\otimes_AM_B)$ and
$\Omega^n_B(T)\otimes_B\check M\otimes_AM_B$ in
$\underline{\mathrm{mod}}(B)$, we obtain that
$\eta_{\Omega^n_B(T)}: \Omega^n_B(T_B)\to
\Omega^n_B(T)\otimes_B\check M\otimes_AM_B$ is an isomorphism  in
$\underline{\mathrm{mod}}(B)$.

As we have an exact sequence of $B$-modules
$$\Omega^n_B(T_B)\stackrel {\eta_{\Omega^n_B(T_B)}}{\to}
\Omega^n_B(T_B)\otimes_B \check M\otimes_A M_B\to
\Omega^n_B(T_B)\otimes_B U_B\to 0,$$ we deduce that
$\eta_{\Omega^n_B(T_B)}$ 
has projective cokernel.
In fact, let $S_B$
be an indecomposable direct summand of $\Omega^n_B(T_B)$.
Then $\eta_{\Omega^n_B(T_B)}$ is
the direct sum of such $\eta_S$ and $\eta_S$ is an isomorphism in  $\underline{\mathrm{mod}}(B)$.
If $S_B$ is projective,  $\eta_S$ is injective  and
$\eta_S$ has projective cokernel,  since
$\eta_B$ is injective with projective cokernel $U_B$ by Claim 1.
If $S_B$ is not projective, then  the fact that  $\eta_S$ is an isomorphism in
$\underline{\mathrm{mod}}(B)$ implies that $\eta_S$
has projective cokernel.

Since $\eta_{\Omega^n_B(T_B)}$
has projective cokernel,
$\Omega^n_B(T_B)\otimes_B U_B$ is projective and  the module ${}_BU_B$ is  strongly  right nonsingular.
\end{Proof}

\begin{Proof}[Proof of Claim 2]
 Let $T_A$ be an $A$-module. For $n\geq 1$,
$\Omega_A^n(T_A)\otimes_AM_B\simeq \Omega_B^n(T\otimes_AM_B)\oplus
P_B$ with $P_B$ projective. Then
$$\Omega_A^n(T)\otimes_AM\otimes_BU\otimes_BN_A\simeq
(\Omega_B^n(T\otimes_AM_B)\otimes_BU\otimes_BN_A)\oplus
(P\otimes_BU\otimes_BN_A).$$ The $A$-module
$\Omega_B^n(T\otimes_AM_B)\otimes_BU\otimes_BN_A$ is projective
for $n$ big enough, as   ${}_BU_B$ is  strongly  right nonsingular
and that ${}_BN_A$ is projective as a  left and right module; the
module $P\otimes_BU\otimes_BN_A$ is projective since $U_B$ is
projective. We have proved that
$\Omega_A^n(T)\otimes_AM\otimes_BU\otimes_BN_A$ is projective for
$n>>0$ and that ${}_AM\otimes_BU\otimes_BN_A$ is  strongly  right
nonsingular.

The fact that ${}_AM\otimes_B\check M\otimes_A X_A$ is
strongly  right nonsingular follows from the fact that ${}_AX_A$ is in
$\Smooth(A^e)$ and that ${}_AM\otimes_B\check M\otimes_A X_A$ is
projective as a left and right module.

\end{Proof}

\begin{Proof}[Proof of Claim 3]
The fact that
${}_BN\otimes_A\tilde{X}_A$ is  strongly  right nonsingular follows
from that ${}_BN_A$ is projective as a  left and right module and
that ${}_A\tilde X_A$ is strongly right  nonsingular.

The fact that  ${}_BY\otimes_B\check M_A$ is   strongly right
nonsingular follows from that ${}_BY_B$ is in $\Smooth(B^e)$
and that ${}_B\check M_A$ is projective as a left and right module.

Suppose $\check M=\check M_1\oplus \check M_2$
  as $B^{op}\otimes_KA$-modules. Then
$\mathrm{Hom}_A({}_B\check M_A, {}_AA_A)\simeq {}_AM_B$
is indecomposable as $B^{op}\otimes_KA$-module
implies that $\mathrm{Hom}_A(\check M_1, {}_AA_A)=0$ or
$\mathrm{Hom}_A(\check M_2, {}_AA_A)=0$. But $M_B$ is projective,
and so this happens only if $M_1=0$ or $M_2=0$.
Therefore ${}_B\check M_A$  is
indecomposable.
\end{Proof}

\medskip

We obtain the analogous result to \cite[Corollary 3.1]{DugasMartinezVilla}.

\begin{Cor}\label{adjoint}
Under the same assumption of Theorem~\ref{AdjointPairs}, suppose
further that $Hom_B({}_AM_B, {}_BB_B)$ is projective as an $A$-module, or
$Hom_{B^{op}}({}_BN_A,  {}_BB_B)$ is projective as a left $A$-module.
 Then
$${}_BN_A\simeq \mathrm{Hom}_{A^{op}}({}_AM_B, {}_AA_A)\simeq \mathrm{Hom}_B({}_AM_B, {}_BB_B)$$
and
$${}_AM_B\simeq
\mathrm{Hom}_{A}({}_BN_A, {}_AA_A) \simeq \mathrm{Hom}_{B^{op}}({}_BN_A, {}_BB_B).$$
Moreover
$(M\otimes_B-,N\otimes_A-)$ and $(N\otimes_A-,M\otimes_B-)$ are
adjoint functors between $\mathrm{mod}(A^{op})$ and
$\mathrm{mod}(B^{op})$, which induce pairs of equivalences of the
corresponding singularity categories.

Finally $(-\otimes_AM_B,-\otimes_BN_A)$ and
$(-\otimes_BN_A, -\otimes_AM_B)$ are adjoint
functors between $\mathrm{mod}(A)$ and $\mathrm{mod}(B)$, which
induce pairs of equivalences of the corresponding singularity
categories.
\end{Cor}

\begin{Proof} As a left (resp. right)
adjoint to a given functor is unique up to isomorphisms,
Theorem~\ref{AdjointPairs} shows that
$\mathrm{Hom}_{A}({}_BN_A, {}_AA_A)\simeq {}_AM_B$ and in particular, $\mathrm{Hom}_{A}({}_BN_A, {}_AA_A)$ is projective as a right $B$-module.

On the other hand, if we suppose in Theorem 3.1 that
$\mathrm{Hom}_{A}({}_BN_A, {}_AA_A)$ is projective as a
right $B$-module instead of being projective for
$\mathrm{Hom}_{A^{op}}({}_AM_B, {}_AA_A)$ as a left $B$-module,
a dual proof as that of  Theorem~\ref{AdjointPairs}, by considering
the functors $(\mathrm{Hom}_{A}({}_BN_A, {}_AA_A)\otimes_B-, N\otimes_A-)$
between left module categories $\mathrm{mod}(B^{op})$ and $\mathrm{mod}(A^{op})$,
gives that
$${}_AM_B\simeq \mathrm{Hom}_{A}({}_BN_A, {}_AA_A)$$
as $A^{op}\otimes_KB$-modules, and $(M\otimes_B-, N\otimes_A-)$ is
a pair of adjoint functors between the module categories $\mathrm{mod}(A^{op})$
and  $\mathrm{mod}(B^{op})$. As in the first paragraph, we see that
$\mathrm{Hom}_{A^{op}}({}_AM_B, {}_AA_A)\simeq {}_BN_A$ and in particular,
$\mathrm{Hom}_{A^{op}}({}_AM_B, {}_AA_A)$ is projective as a left $B$-module.

This shows that the the condition that
$\mathrm{Hom}_{A^{op}}({}_AM_B, {}_AA_A)$ is projective as a
left $B$-module and the condition that $\mathrm{Hom}_{A}({}_BN_A, {}_AA_A)$
is projective as a right $B$-module  are equivalent. Furthermore, under these
two equivalent conditions, we know that
\begin{itemize}
\item[(i)]  ${}_BN_A\simeq \mathrm{Hom}_{A^{op}}({}_AM_B, {}_AA_A)$
and ${}_AM_B\simeq \mathrm{Hom}_{A}({}_BN_A, {}_AA_A).$

\item[(ii)]
$(M\otimes_B-,N\otimes_A-)$ is a pair of
adjoint functors  between $\mathrm{mod}(B^{op})$ and
$\mathrm{mod}(A^{op})$, which induce pairs of equivalences of the
corresponding singularity categories.

\item[(iii)]
$(-\otimes_BN_A, -\otimes_AM_B)$ is a pair of  adjoint
functors between $\mathrm{mod}(B)$ and $\mathrm{mod}(A)$, which
induce pairs of equivalences of the corresponding singularity
categories.
\end{itemize}

A dual proof of the above argument shows that  the condition that
$Hom_B({}_AM_B, {}_BB_B)$ is projective as an $A$-module, and the condition that
$Hom_{B^{op}}({}_BN_A,  {}_BB_B)$ is projective as an $A$-module, are
equivalent; under these two conditions, we have
\begin{itemize}
\item[(i)]  ${}_BN_A\simeq \mathrm{Hom}_{B}({}_AM_B, {}_BB_B)$
and
${}_AM_B\simeq \mathrm{Hom}_{B^{op}}({}_BN_A, {}_BB_B).$

\item[(ii)]
$({}_BN\otimes_A-, {}_AM\otimes_B-)$ is a pair of
adjoint functors  between $\mathrm{mod}(A^{op})$ and
$\mathrm{mod}(B^{op})$, which induce pairs of equivalences of the
corresponding singularity categories.

\item[(iii)]
$(-\otimes_AM_B, -\otimes_BN_A)$ is a pair of  adjoint
functors between $\mathrm{mod}(A)$ and $\mathrm{mod}(B)$, which
induce pairs of equivalences of the corresponding singularity
categories.
\end{itemize}

\end{Proof}

Let $\nu_A:=Hom_K(Hom_A(-,A),K)$ be the Nakayama functor on
$\mathrm{mod}(A)$. If $Q_A$ is a projective $A$-module, then
$\nu(Q_A)$ is an injective $A$-module, and if $I_A$ is an
injective $A$-module, then $\nu(I_A)$ is a projective $A$-module.

\begin{Lem}\label{nakayamacommuteswithsingularequiv}
Under the same assumption of Theorem~\ref{AdjointPairs},  if $I$
is injective as a $B$-module, then $M\otimes_BI$ is injective as an
$A$-module. Moreover $(M\otimes_B-)\circ\nu_B\simeq\nu_A\circ
(M\otimes_B-)$.
\end{Lem}

\begin{Proof} We know that $N\simeq Hom_A(M,A)$ and that $N\otimes_A-$ is
(left and) right adjoint to $M\otimes_B-$. Hence for an injective
$A$-module $I$ we get
$$Hom_{B}(-,N\otimes_AI)\simeq Hom_{A}(M\otimes_B-,I)$$
by Corollary~\ref{adjoint}. Moreover $M\otimes_B-$ is exact since
$M$ is projective as a $B$-module. $Hom_{A}(-,I)$ is exact since $I$
is injective as an $A$-module. Therefore $Hom_{B}(-,N\otimes_AI)$ is
exact as a functor $B-mod\lra (A-mod)^{op}$, and we get therefore
that $N\otimes_AI$ is injective.

We have
\begin{eqnarray*}
Hom_A(M\otimes_B-,A)&\simeq &Hom_B(-,Hom_A(M,A))\\
&\simeq &Hom_B(-,N)\\
&\simeq&Hom_B(-,B)\otimes_BN
\end{eqnarray*}
as right $A$ modules,
since $N$ is projective as $B$-module.
Hence,
\begin{eqnarray*}
\nu_A(M\otimes_B-)&=&Hom_K(Hom_A(M\otimes_B-,A),K)\\
&\simeq &Hom_K(Hom_B(-,B)\otimes_BN,K)\\
&\simeq&Hom_B(N,Hom_K(Hom_B(-,B),K))\\
&\simeq&Hom_B(N,B)\otimes_BHom_K(Hom_B(-,B),K)\\
&\simeq&M\otimes_B\nu_B(-)
\end{eqnarray*}
This shows the lemma. \end{Proof}

\begin{Cor}
Under the same assumption of Theorem~\ref{AdjointPairs},  the functor
$-\otimes_AM_B$ sends projective injective $A$-modules to
projective injective $B$-modules.
\end{Cor}

\section{Singular equivalences of Morita type and Hochschild homology}

\label{hochschildsection}

In this section, we consider invariant property of Hochschild homology under singular equivalences of Morita type.
For stable equivalences of Morita type, in \cite{LiuXi2005}, Yu-Ming Liu
and Chang-Chang Xi proved that a stable equivalence of Morita type preserves
Hochschild homology groups of positive degrees.
Remark that by \cite[Theorem 1.1]{LZZ} the invariance of degree zero Hochschild
homology group under a stable equivalence of Morita type is equivalent to the
famous Auslander-Reiten conjecture on the invariance of the number of non
projective simple modules under
stable equivalence.

We shall now prove that  a singular equivalence  of Morita type
induces an isomorphism of Hochschild homology in positive degrees.

\begin{Thm}\label{Hochschildhomology}
Let $K$ be a Noetherian commutative ring and let
$A$ and $B$ be Noetherian
$K$-algebras which are projective as $K$-modules.
Suppose that $({}_AM_B,{}_BN_A)$ induce a singular
equivalence of Morita type.
\begin{enumerate}
\item
Then there is $n_0\in\N$
so that $HH_n(A)\simeq HH_n(B)$
for each $n>n_0$.
\item
If $K$ is a field, and if $A$ and $B$ are finite dimensional, then
$HH_n(A)\simeq HH_n(B)$
for each $n>0$.
\end{enumerate}
\end{Thm}

Our proof of the first statement, inspired by \cite[Section 1.2]{Zimmermann2007},
is similar to that of \cite[Theorem 4.4]{LiuXi2005}, which uses a change-of-rings
argument. Notice that our argument is simpler than the proof in \cite{LiuXi2005}
and in fact works also
for stable equivalences of Morita type.   Our proof of the second statement
makes use of transfer maps  and is similar to that of \cite[Remark 3.3]{LZZ}.

\begin{Proof}[Proof of Theorem~\ref{Hochschildhomology}.(1).]
Let ${\mathbb B}A$ be the bar resolution of $A$, that is
$${\mathbb B}A:\;\;\dots\lra A^{\otimes 5}\lra A^{\otimes 4}\lra
A^{\otimes 3}\lra A^{\otimes 2}(\lra A\lra 0).$$
Then, we my apply $N\otimes_A-\otimes_AM$ and obtain an exact sequence
$N\otimes_A{\mathbb B}A\otimes_AM$ of $B^e$-modules:
$$\;\;\dots\lra
N\otimes_AA^{\otimes 4}\otimes_AM\lra N\otimes_AA^{\otimes
3}\otimes_AM\lra N\otimes_AA^{\otimes 2}\otimes_AM(\lra
N\otimes_AM\lra 0).$$ of $B^e$-modules, since $M$ and $N$ are
projective on the right, resp. on the left.

The key observation is the following isomorphism of complexes
$$ \begin{array}{rcl}  ({}_AM\otimes_B N_A)\otimes_{A^e}
\mathbb{B}A& \simeq & B\otimes_{B^e} ({}_BN\otimes_A
\mathbb{B}A\otimes_AM_B)\\
(m\otimes n)\otimes u&\mapsto & 1\otimes (n\otimes u\otimes
m).\end{array}$$ which is easily verified. Taking homology groups
gives $$HH_n(A)\oplus Tor_n^{A^e}(X, A)\simeq
Tor_n^{A^e}(M\otimes_B N, A)\simeq  Tor_n^{B^e}(B,
N\otimes_A M)\simeq  HH_n(B)\oplus Tor_n^{B^e}(B, Y)$$ for each
$n\geq 0$. When $n$ is large, $Tor_n^{A^e}(X, A)\simeq 0\simeq
Tor_n^{B^e}(B, Y)$, as $X\in \Smooth(A^e)$ and $Y\in
\Smooth(B^e)$, we obtain that $HH_n(A)\simeq HH_n(B)$ for $n>>0$.
\end{Proof}

For the proof of Theorem~\ref{Hochschildhomology}.(2),
let us recall some properties of transfer
maps in Hochschild homology.

Let $A$ and $B$ be two algebras over a commutative ring $k$. Let
$M$ be an $A$-$B$-bimodule such that $M_B$ is finitely generated
and projective. Then we can define  a transfer map $t_M:
HH_n(A)\to HH_n(B)$ for each $n\geq 0$. As we don't need the
construction of this map, we refer the reader to Bouc~\cite{Bouc}
(see also \cite{LZZ,KLZ} for a summary of Bouc's results).

\begin{Prop} \label{Bouc}\cite[Section 3]{Bouc}
Let $A$, $B$ and $C$ be $k$-algebras over a commutative ring $k$.
\begin{enumerate}
\item If $M$ is an $A$-$B$-bimodule and $N$ is a $B$-$C$-bimodule
such that $M_B$ and $N_C$ are finitely generated and projective,
then we have $t_N\circ t_M=t_{M\otimes_B N}: HH_n(A)\rightarrow
HH_n(C)$, for each $n\geq 0$.

\item Let $$0\rightarrow L\rightarrow M\rightarrow N\rightarrow
0$$ be a short exact sequence of $A$-$B$-bimodules which are
finitely generated and projective as right $B$-modules. Then
$t_M=t_L+t_N: HH_n(A)\rightarrow HH_n(B)$, for each $n\geq 0$.

\item Suppose that $k$ is an algebraically closed field and that
$A$ and $B$ are finite dimensional $k$-algebras. Then for a
finitely generated projective $A$-$B$-bimodule $P$, the transfer
map $t_P: HH_n(A)\rightarrow HH_n(B)$ is zero for each $n>0$.

\item Consider $A$ as an $A$-$A$-bimodule by left and right
multiplications, then $t_A: HH_n(A)\rightarrow HH_n(A)$ is the
identity map for any $n\geq 0$.
\end{enumerate}
\end{Prop}

\begin{Proof}[Proof of Theorem~\ref{Hochschildhomology}.(2).]
For $n\geq 0$, we have transfer maps
$t_M: HH_n(A)\to HH_n(B)$ and  $t_N: HH_n(B)\to HH_n(A)$. By the
above result, $$t_N\circ t_M=t_{M\otimes_BN}=t_A+t_X=Id+t_X$$ as
maps from $HH_n(A)$ to itself.

Let $\overline{K}$ be the algebraic closure of $K$ and write
\begin{eqnarray*}
\overline{A}&=&A\otimes_K\overline{K},\\
\overline{B}&=&B\otimes_K\overline{K},\\
\overline{M}&=&M\otimes_K\overline{K},\\
\overline{N}&=&N\otimes_K\overline{K},\\
\overline{X}&=&X\otimes_K\overline{K},\\
\overline{Y}&=&Y\otimes_K\overline{K}.
\end{eqnarray*}
 Then one verifies easily
that $({}_{\ol{A}}\ol{M}_{\ol{B}}, {}_{\ol{B}}\ol{N}_{\ol{A}})$
induces a singular equivalence of Morita type between $\ol{A}$ and
$\ol{B}$, because
$${}_{\ol{A}}\ol{M}\otimes_{\ol{B}}\ol{N}_{\ol{A}}\simeq
{}_{\ol{A}}\ol{A}_{\ol{A}}\oplus {}_{\ol{A}}\ol{X}_{\ol{A}}$$ with
 $\ol{X}\in \mathcal{P}^{<\infty}(\ol{A}^e)$;
$${}_{\ol{B}}\ol{N}\otimes_{\ol{A}}\ol{M}_{\ol{B}}\simeq
{}_{\ol{B}}\ol{B}_{\ol{B}}\oplus {}_{\ol{B}}\ol{Y}_{\ol{B}}$$ with
$\ol{Y}\in \mathcal{P}^{<\infty}(\ol{B}^e)$. We also have
$t_{\ol{M}}=t_M\otimes_Kid_{\ol{K}}$.

Since $\ol{X}\in \mathcal{P}^{<\infty}(\ol{A}^e)$, there is an
exact sequence of $\ol A^e$-modules
$$0\to \ol P_n\to \cdots \to \ol P_0\to \ol{X}\to 0$$
with $\ol P_0, \cdots, \ol P_n$ projective.
By the point (2)(3) of Proposition~\ref{Bouc},  for
$n>0$, we have $t_{\ol{X}}=\sum_{i=0}^n(-1)^it_{\ol P_i}=0$ as a
homomorphism from $ HH_n(\ol{A})\to HH_n(\ol{A})$,   and thus
$t_X=0: HH_n(A)\to HH_n(A)$ for $n>0$. This shows that
$$t_N\circ t_M: HH_n(A)\to
HH_n(A)$$ and $$t_M\circ t_N: HH_N(B)\to
HH_n(B)$$ are isomorphisms for $n>0$. We deduce that $$t_M: HH_n(A)\to
HH_n(B)$$ is an isomorphism for $n>0$.
\end{Proof}

\begin{Rem} 
Finally we briefly mention what is known in this context about
invariance of Hochschild cohomology under stable equivalence
of Morita type and under singular equivalence of Morita type.

Chang-Chang Xi
prove in \cite[Theorem 4.2]{Xi2008} that a stable equivalence
of Morita type between Artin algebras
preserves the Hochschild cohomology groups of positive degrees,
generalising a previous result of Zygmunt Pogorza{\l}y
\cite[Theorem 1.1]{Pogorzaly} for selfinjective algebras.
Sheng-Yong Pan and the first author further showed in
\cite{PanZhou} the invariance of
stable Hochschild cohomology rings under stable equivalences of Morita type.

Chen and Sun
prove in \cite{ChenSun} that  Tate-Hochschild cohomology
rings of Gorenstein algebras are preserved under singular equivalences of Morita type.
A careful study of the proof of  \cite[Theorem 4.2]{Xi2008}
shows that the proof  of \cite[Theorem 4.2]{Xi2008} works for
singular equivalences of Morita type. We obtain from this study
that a singular equivalence of Morita type preserves
Hochschild cohomology groups of large degrees. However, we
don't know the algebra structure.
\end{Rem}

\end{document}